\documentclass[12pt]{article}
\usepackage{amsmath,amssymb,amsfonts,epsfig}
\author{S.~V.~Chmutov\thanks{The Ohio State University, Mansfield.}, S.~K.~Lando
\thanks{Institute for System Research RAS and the Poncelet Laboratory, Independent
University of Moscow, partly supported by the grant
ACI-NIM-2004-243 (Noeuds et tresses), RFBR 05-01-01012-a, NWO-RFBR
047.011.2004.026 (RFBR 05-02-89000-NWOa), GIMP
ANR-05-BLAN-0029-01.}}

\title{Mutant knots and intersection graphs}
\date{April 10, 2007}

\usepackage{theorem}

\def\C{{\mathbb C}}

\def\Z{{\mathbb Z}}

\def\cF{{\cal F}}

\def\cG{{\cal G}}

\def\cP{{\cal P}}

\newcommand{\sL}{\mathfrak{sl}} 
\newcommand{\gl}{\mathfrak{gl}}

\newtheorem{theorem}{Theorem}
\newtheorem{lemma}{Lemma}

{\theorembodyfont{\rmfamily}    

\newtheorem{definition}{Definition}

}

\newcommand{\rb}{\raisebox}
\newcommand{\ig}{\includegraphics}
\newcommand\risS[6]{\rb{#1pt}[#5pt][#6pt]{\begin{picture}(#4,15)(0,0)
  \put(0,0){\ig[width=#4pt]{#2.eps}} #3
     \end{picture}}}

\begin{document}

\maketitle

\begin{abstract}
\noindent We prove that if a finite order knot invariant does not
distinguish mutant knots, then the corresponding weight system
depends on the intersection graph of a chord diagram rather than
on the diagram itself. The converse statement is easy and well
known. We discuss relationship between our results and certain Lie
algebra weight systems.
\end{abstract}

\section{Introduction}

Below, we use standard notions of the theory of finite order, or
Vassiliev, invariants of knots in $3$-space; their definitions can
be found, for example, in~\cite{CDbook} or~\cite{LZ}. All knots
are assumed to be oriented.

Two knots are said to be {\it mutant} if they differ by a rotation/reflection
of a tangle with four endpoints; if necessary, the orientation
inside the tangle may be replaced by the opposite one. Here is a
famous example of mutant knots, the Conway ($11n34$) knot $C$ of
genus 3, and Kinoshita--Terasaka ($11n42$) knot $KT$ of genus 2
(see \cite{Atlas}).\label{con-knot}
$$C =\ \risS{-20}{c-kn}{}{60}{20}{20}
\hspace{2cm}
  KT =\ \risS{-20}{tk-kn}{}{60}{20}{20}
$$
Note that the change of the orientation of a knot can be achieved
by a mutation in the complement to a trivial tangle.

Most known knot invariants cannot distinguish mutant knots.
Neither the (colored) Jones polynomial, nor the HOMFLY polynomial,
nor the Kauffman two variable polynomial distinguish mutants.  All
Vassiliev invariants up to order 10 do not distinguish mutants as
well \cite{Mu} (up to order 8 this fact was established by a
direct computation \cite{CDL1,CDbook}). However, there is a
Vassiliev invariant of order 11 distinguishing $C$ and $KT$
\cite{MC,Mu}. It comes from the colored HOMFLY polynomial.

The main combinatorial objects of the Vassiliev theory of knot
invariants are {\it chord diagrams}. To a chord diagram, its {\it
intersection graph} (also called {\it circle graph}) is
associated. The vertices of the graph correspond to chords of the
diagram, and two vertices are connected by an edge if and only if
the corresponding chords intersect.

The value of a Vassiliev invariant of order $n$ on a singular knot
with $n$ double points depends only on the chord diagram of the
singular knot. Hence any such invariant determines a function, a
{\it weight system}, on chord diagrams with $n$ chords.
Conversely, any weight system induces, in composition with the
Kontsevich integral, which is the universal finite order
invariant, a finite order invariant of knots. Such knot invariants
are called {\it canonical}. Canonical invariants span the whole
space of Vassiliev invariants.

Direct calculations for small $n$
show that the values of these functions are uniquely determined by the intersection graphs of the chord diagrams. This fact motivated the {\it intersection graph conjecture} in \cite{CDL1} (see also \cite{CDbook}) which states that any weight system depends on the intersection graph only. This conjecture happened to be false, because of the
existence of a finite order invariant that distinguishes two mutant knots
mentioned above and the following fact.

{\it The knot invariant induced by a weight system whose values depend only on the intersection graph of the chord diagrams cannot distinguish mutants.}

A justification of this statement, due to T.~Le (unpublished), looks like follows
(see details in~\cite{CDbook}). If we have a knot (in general position)
with a distinguished
two-string tangle, then all the terms in the Kontsevich integral of the knot
having chords connecting the tangle with its exterior vanish.

Our goal is to prove the converse statement thus establishing an equivalence
between finite order knot invariants nondistinguishing mutants and
weight systems depending on the intersection graphs of chord diagrams only.

\begin{theorem}\label{main-theorem}
If a finite order knot invariant does not distinguish mutants, then
the corresponding weight system depends only on the intersection graphs
of chord diagrams.
\end{theorem}

Together, the two statements can be combined as follows.

{\it A canonical knot invariant does not distinguish mutants if
and only if its weight system depends on the intersection graphs
of chord diagrams only.}

Recently, B.~Mellor \cite{Mel} extended the concept of
intersection graph to string links. We do not know whether our
Theorem \ref{main-theorem} admits an appropriate generalization.

Section~\ref{s:proof} is devoted to the proof of Theorem~1. In Sec.~\ref{s:Lie}, we
discuss relationship between intersection graphs and the weight
systems associated to the Lie algebra~$\sL(2)$ and the Lie
algebra~$\gl
(1|1)$.

The paper was written during the second author's
visit to the Mathematical Department 
of the Ohio State University. He expresses
his gratitude to this institution for warm hospitality
and excellent working conditions.
The authors are grateful to S.~Duzhin, K.~J.~Supowit, and A.~Vaintrob
for useful discussions.

\section{Proof}\label{s:proof}

\subsection{Representability of graphs as the intersection graphs
of chord diagrams}

Not every graph can be represented as the intersection graph of a chord diagram. For example, the following graphs are not intersection graphs.
$$\ig[width=15mm]{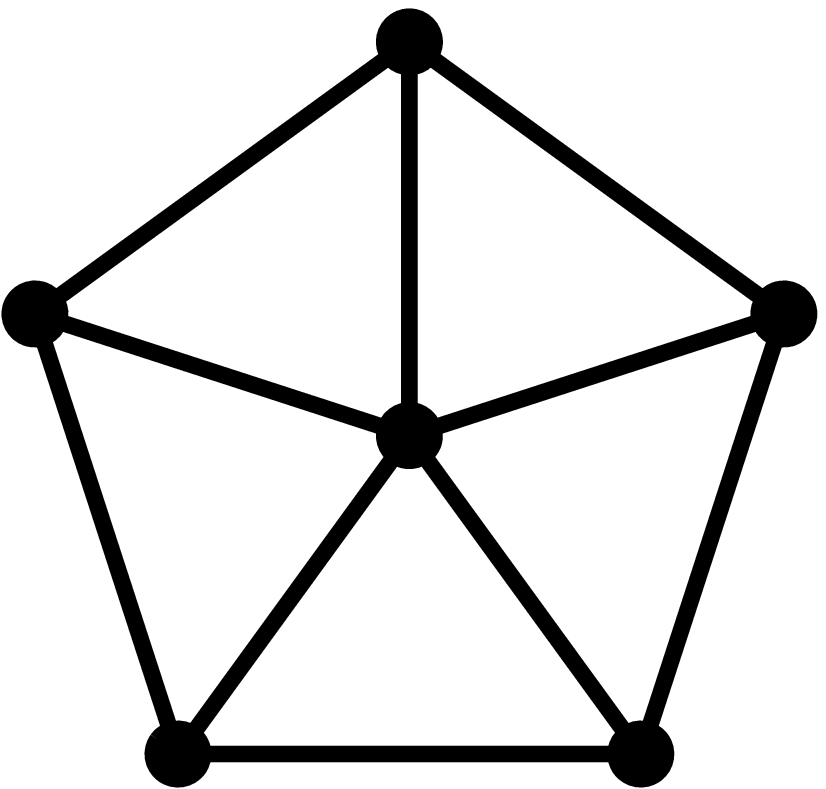}\hspace{2cm}
  \ig[width=13.5mm]{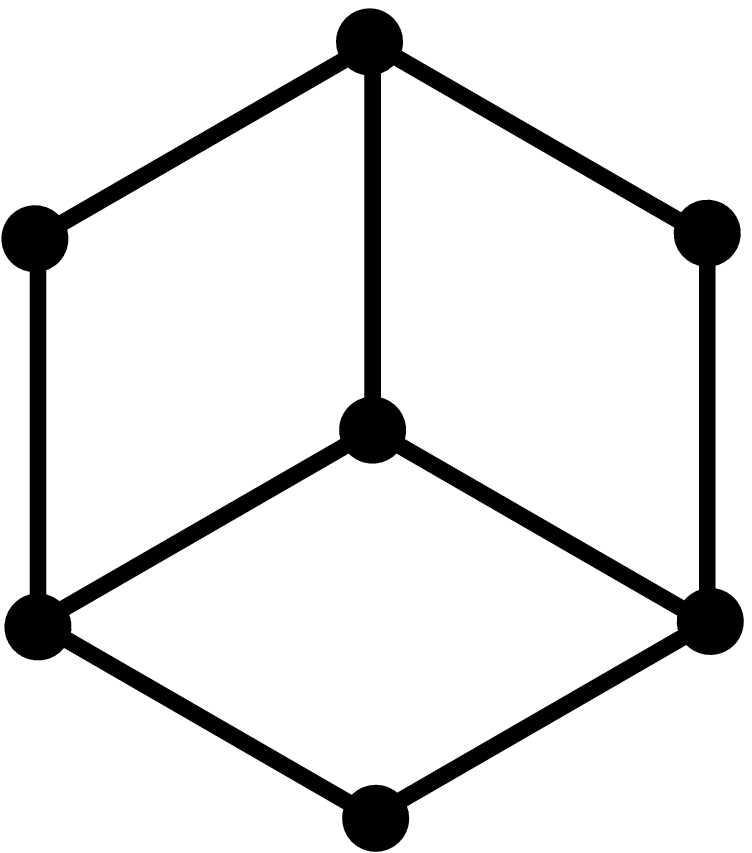}\hspace{2cm}
  \ig[width=15mm]{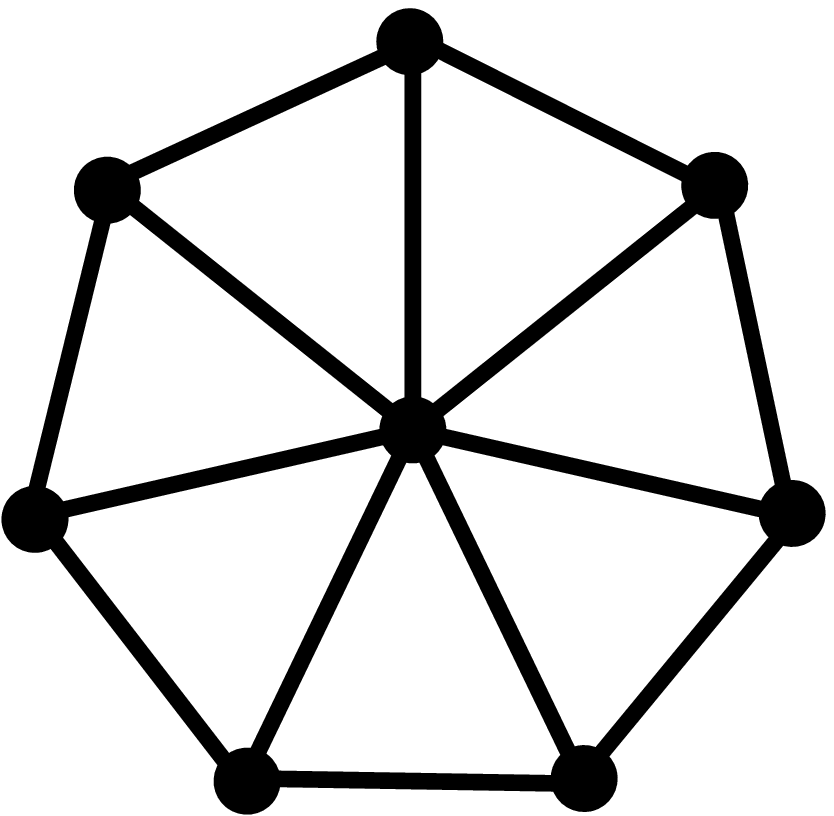}
$$
A characterization of those graphs that can be realized as
intersection graphs is given by an elegant theorem of A.~Bouchet
\cite{B2}.

On the other hand, distinct diagrams may have coinciding intersection graphs.
For example, next three diagrams have the same intersection graph
$\ig[width=25mm]{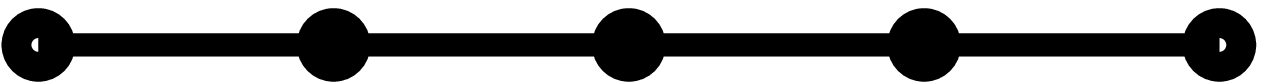}$\ :
$$  \ig[width=15mm]{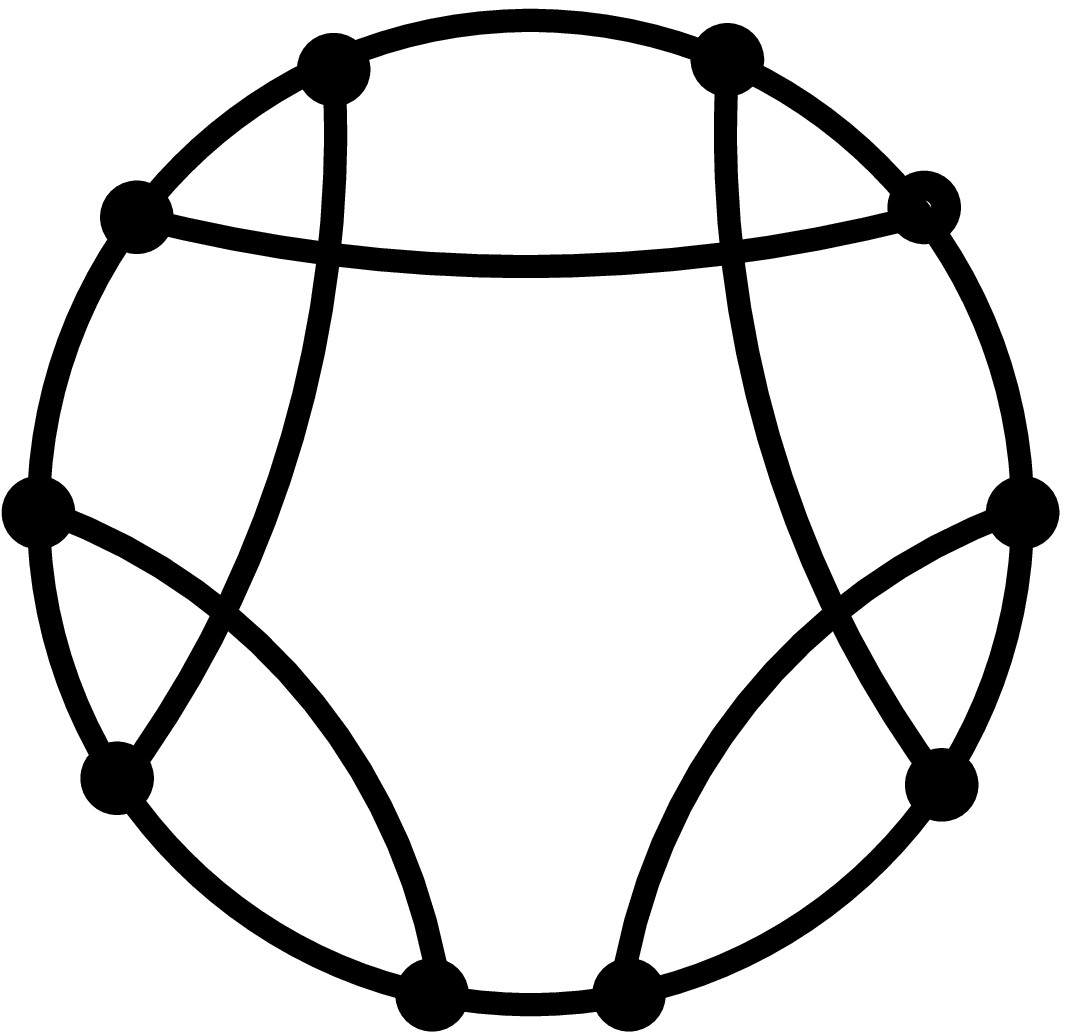}\qquad
  \ig[width=15mm]{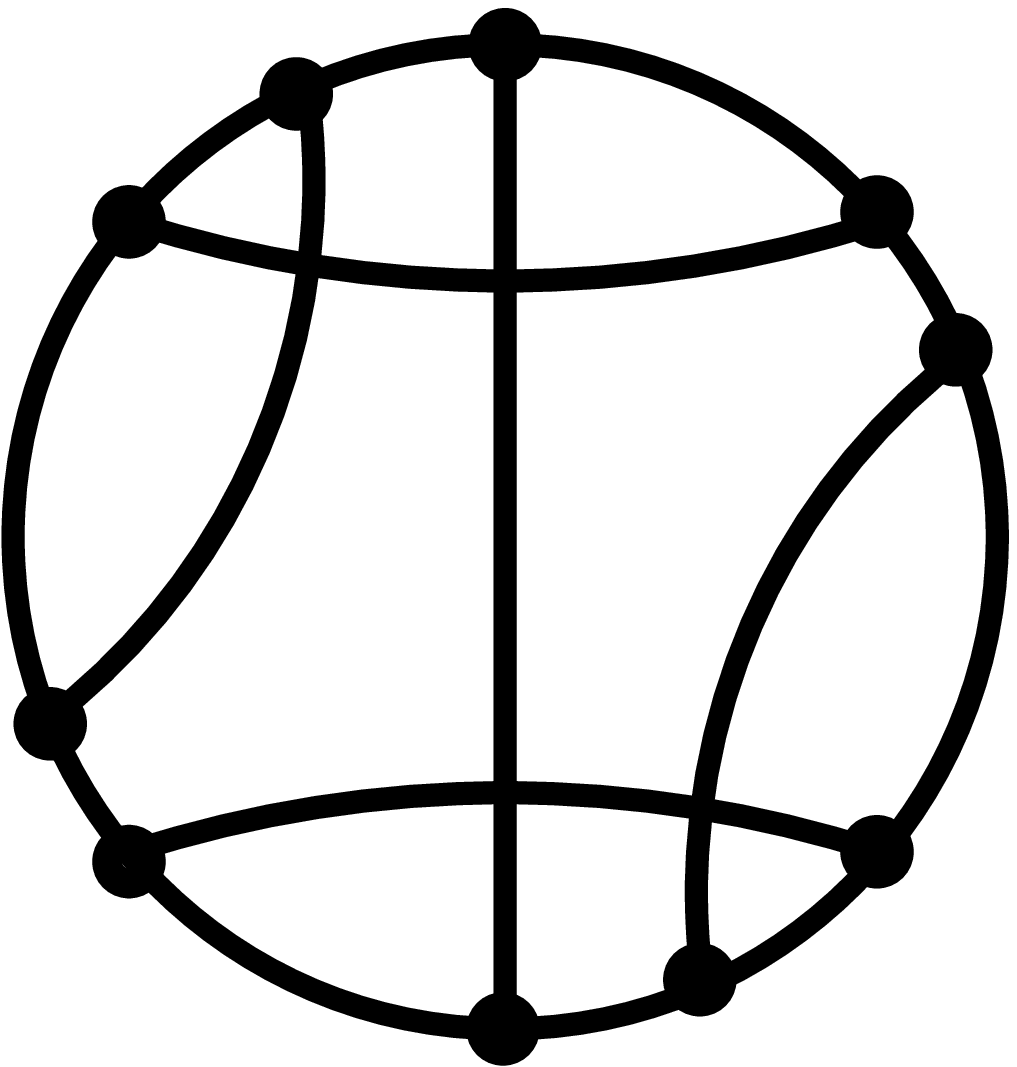}\qquad
  \ig[width=15mm]{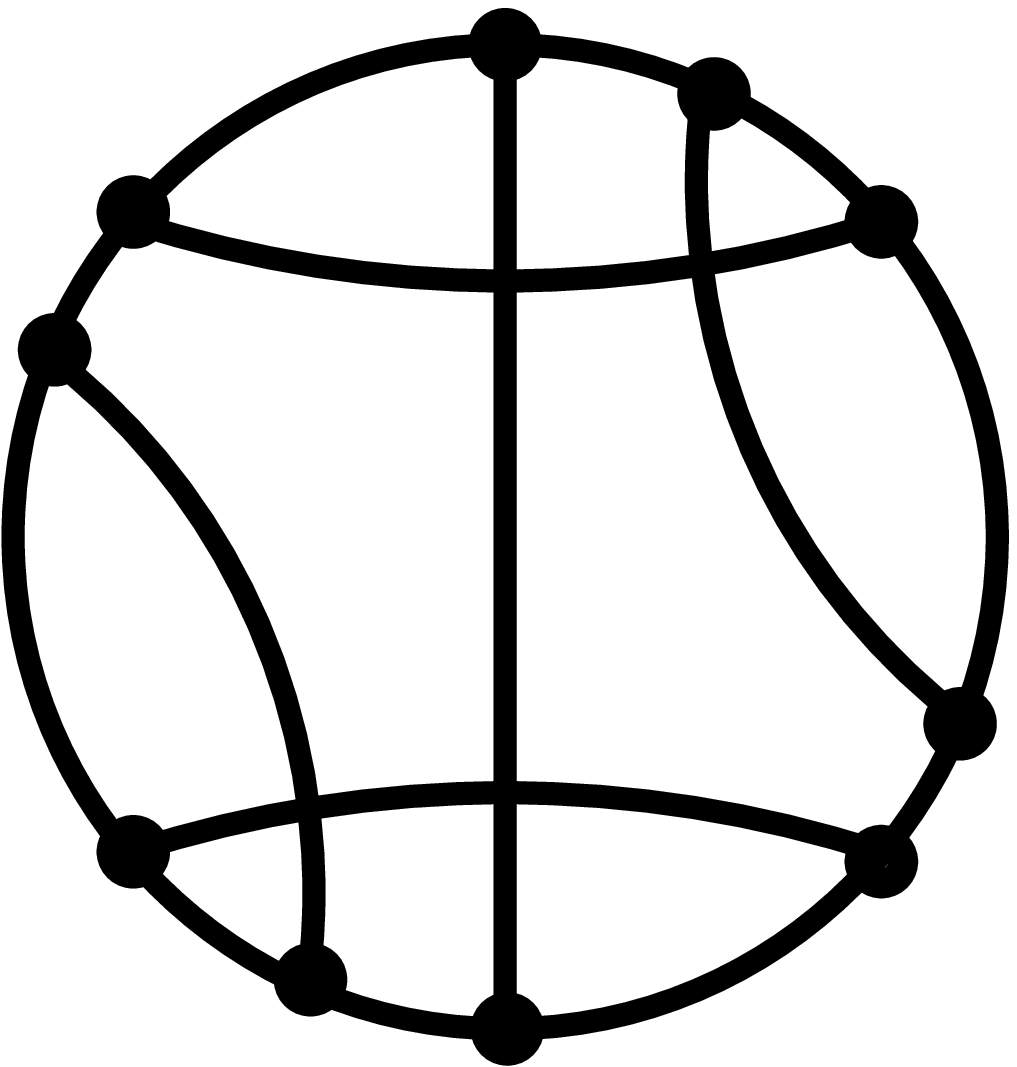}
$$

A combinatorial analog of the tangle in
mutant knots is a {\it share} \cite{CDL1,CDbook}.
Informally, a {\it share} of a chord diagram is a subset
of chords whose endpoints are separated into at most two parts by the endpoints of the complementary chords. More formally,
\begin{definition}
A \textit{share} is a part of a chord diagram consisting
of two arcs of the outer circle possessing the following property:
each chord one of whose ends belongs to these arcs
has both ends on these arcs.
\end{definition}

Here are some examples:
$$\risS{-20}{share1}{
       \put(5,-10){\mbox{\scriptsize A share}}}{40}{25}{30}\hspace{2cm}
  \risS{-20}{share2}{
       \put(-3,-10){\mbox{\scriptsize Not a share}}}{40}{20}{20}\hspace{2cm}
  \risS{-20}{share3}{
       \put(0,-10){\mbox{\scriptsize Two shares}}}{48}{20}{20}
$$
The complement of a share also is a share.
The whole chord diagram is its own share whose complement contains no chords.

\begin{definition}
A \textit{mutation of a chord diagram} is another chord diagram obtained by a rotation/reflection of a share.
\end{definition}

For example, three mutations of the share in the first chord diagram above produce
the following  chord diagrams:\label{st-example}
$$\risS{-20}{sh-mut1}{}{40}{25}{20}\hspace{2cm}
  \risS{-20}{sh-mut2}{}{40}{20}{20}\hspace{2cm}
  \risS{-20}{sh-mut3}{}{40}{20}{20}
$$

Obviously, mutations preserve the intersection graphs of chord diagrams.

\begin{theorem} \label{cd-mut-theorem}
Two chord diagrams have the same intersection graph if and only if
they are related by a sequence of mutations.
\end{theorem}

This theorem is contained implicitly in papers \cite{B1,Co,GSH}
where chord diagrams are written as {\it double occurrence words},
the language better suitable for describing algorithms than for
topological explanation.

\bigskip
{\bf Proof of Theorem \ref{cd-mut-theorem}.}\\
The proof of this theorem uses Cunningham's theory of graph decompositions \cite{Cu}.

A {\it split} of a (simple) graph $\Gamma$ is a disjoint
bipartition $\{V_1,V_2\}$ of its set of vertices $V(\Gamma)$ such
that each part contains at least 2 vertices, and there are subsets
$W_1\subseteq V_1$, $W_2\subseteq V_2$ such that all the edges of
$\Gamma$ connecting $V_1$ with $V_2$ form the complete bipartite
graph $K(W_1,W_2)$ with the parts $W_1$ and $W_2$. Thus for a
split $\{V_1,V_2\}$ the whole graph $\Gamma$ can be represented as
a union of the induced subgraphs $\Gamma(V_1)$ and $\Gamma(V_2)$
linked by a complete bipartite graph.

Another way to think about splits, which is sometimes more
convenient and which we shall use in the pictures below, looks
like follows. Consider two graphs $\Gamma_1$ and $\Gamma_2$ each
having a distinguished vertex $v_1\in V(\Gamma_1)$ and $v_2\in
V(\Gamma_2)$, respectively, called {\it markers}. Construct the
new graph $\Gamma = \Gamma_1 \boxtimes_{(v_1,v_2)} \Gamma_2$ whose
set of vertices is
$$V(\Gamma) = \{V(\Gamma_1)-v_1\}\sqcup \{V(\Gamma_2)-v_2\}$$
and whose set of edges is
$$\begin{array}{rcl}
E(\Gamma) &=&
    \{(v'_1,v''_1)\in E(\Gamma_1) : v'_1\not=v_1\not=v''_1 \}\ \sqcup\
    \{(v'_2,v''_2)\in E(\Gamma_2) : v'_2\not=v_2\not=v''_2 \}\ \sqcup \vspace{8pt}\\
&&  \{(v'_1,v'_2): (v'_1,v_1)\in E(\Gamma_1)\ \mbox{and}\ (v_2,v'_2)\in E(\Gamma_2)\}\ .
\end{array}$$
Representation of $\Gamma$ as $\Gamma_1 \boxtimes_{(v_1,v_2)}
\Gamma_2$ is called a {\it decomposition} of $\Gamma$, $\Gamma_1$
and $\Gamma_2$ are called the {\it components} of the
decomposition. The partition $\{V(\Gamma_1)-v_1,
V(\Gamma_2)-v_2\}$ is a split of $\Gamma$. Graphs $\Gamma_1$ and
$\Gamma_2$ might be decomposed further giving a finer
decomposition of the initial graph $\Gamma$. Pictorially, we
represent a decomposition by pictures of its components where the
corresponding markers are connected by a dashed edge.

A {\it prime} graph is a graph with at least three vertices
admitting no splits. A decomposition of a graph is said to be {\it
canonical} if the following conditions are satisfied:
\begin{itemize}
\item[(i)] each component is either a prime graph, or a complete
graph $K_n$, or a star $S_n$, which is the tree with a vertex, the
{\it center}, adjacent to $n$ other vertices; \item[(ii)] no two
components that are complete graphs are neighbors, that is, their
markers are not connected by a dashed edge; \item[(iii)] the
markers of two components that are star graphs connected by a
dashed edge are either both centers or both not centers of their
components.
\end{itemize}

W.~H.~Cunningham proved \cite[Theorem 3]{Cu} that each graph with
at least six vertices possesses a unique canonical decomposition.

Let us illustrate the notions introduced above by two examples of
canonical decomposition of the intersection graphs of chord
diagrams. We number the chords and the corresponding vertices in
our graphs, so that the unnumbered vertices are the markers of the
components. The first example is our example from
page~\pageref{st-example}:
$$\risS{-20}{cd-st-ex}{
       \put(0,-10){\mbox{\scriptsize A chord diagram}}}{55}{30}{40}\hspace{2cm}
  \risS{-15}{ig-st-ex}{
       \put(-3,-15){\mbox{\scriptsize The intersection graph}}}{60}{20}{20}\hspace{2cm}
  \risS{-20}{candec-st-ex}{
       \put(20,-10){\mbox{\scriptsize The canonical decomposition}}}{140}{20}{20}
$$
The second example represents the chord diagram of the double
points in the plane diagram of the Conway knot $C$ from page
\pageref{con-knot}. The double points of the shaded tangle are
represented by the chords 1,2,9,10,11.
$$
\risS{-25}{ccd}{
       \put(0,-10){\mbox{\scriptsize Chord diagram}}}{55}{30}{45}\hspace{1.5cm}
  \risS{-25}{ig-c-ex}{
       \put(20,-10){\mbox{\scriptsize Intersection graph}}}{100}{20}{20}\hspace{1.5cm}
  \risS{-23}{candec-c-ex}{
       \put(60,-12){\mbox{\scriptsize Canonical decomposition}}}{180}{20}{20}
$$

The key observation in the proof of Theorem \ref{cd-mut-theorem}
is that components of the canonical decomposition of any
intersection graph admit a unique representation by chord
diagrams. For a complete graph and star components, this is
obvious. For a prime component, this was proved by A.~Bouchet
\cite[Statement 4.4]{B1} (see also \cite[Section 6]{GSH} for an
algorithm finding such a representation for a prime graph).

Now to describe all chord diagrams with a given intersection
graph, we start with a component of its canonical decomposition.
There is only one way to realize the component by a chord diagram.
We draw the chord corresponding to the marker as a dashed chord
and call it the {\it marked chord}. This chord indicates the
places where we must cut the circle removing the marked chord
together with small arcs containing its endpoints. As a result we
obtain a chord diagram on two arcs. Repeating the same procedure
with a neighbor component of the canonical decomposition, we get
another chord diagram on two arcs. We have to sew these two
diagrams together by their arcs in an alternating order. There are
four possibilities to do this, and they differ by mutations of the
share corresponding to the second (or, alternatively, the first)
component. This completes the proof of Theorem
\ref{cd-mut-theorem}. \hspace{\fill}$\square$

\medskip
To illustrate the last stage of the proof consider our standard example and take the star
2-3-4 component first and then the triangle component. We get
$$\risS{-12}{re-st-com}{}{160}{25}{8} \qquad\mbox{and}\qquad
  \risS{-12}{re-tr-com}{}{160}{25}{8}
$$
Because of the symmetry, the four ways of sewing these diagrams
produce only two distinct chord diagrams with a marked chord:
$$\risS{-12}{re-st-tr}{}{100}{15}{8} \qquad\mbox{and}\qquad
  \risS{-12}{re-st-rt}{}{100}{15}{8};
$$
repeating the same procedure with the marked chord for the last
1-6 component of the canonical decomposition, we get
$$\risS{-12}{re-os-com}{}{160}{20}{10}
$$
Sewing this diagram into the previous two in all possible ways we
get four mutant chord diagrams from page \pageref{st-example}.

As an enjoyable exercise we leave to the reader to work out our
second example with the chord diagram of the diagram of the Conway
knot and find the mutation producing the chord diagram of the
plane diagram of the Kinoshita--Terasaka knot using the canonical
decomposition.

\subsection{Proof of Theorem \ref{main-theorem}}

Suppose we have a Vassiliev knot invariant $v$ of order at most
$n$ that does not distinguish mutant knots. Let $D_1$ and $D_2$ be
chord diagrams with $n$ chords whose intersection graphs coincide.
We are going to prove that the values of the weight system of $v$
on $D_1$ and $D_2$ are equal.

By Theorem \ref{cd-mut-theorem}, it is enough to consider the case
when $D_1$ and $D_2$ differ by a single mutation in a share $S$.
Let $K_1$ be a singular knot with $n$ double points whose chord
diagram is $D_1$. Consider the collection of double points of
$K_1$ corresponding to the chords occurring in the share $S$. By
the definition of a share, $K_1$ has two arcs containing all these
double points and no others. By sliding the double points along
one of these arcs and shrinking the other arc we may enclose these
arcs into a ball whose interior does not intersect the rest of the
knot. In other words, we may isotope the knot $K_1$ to a singular
knot so as to collect all the double points corresponding to $S$
in a tangle $T_S$. Performing an appropriate rotation of $T_S$ we
obtain a singular knot $K_2$ with the chord diagram $D_2$. Since
$v$ does not distinguish mutants, its values on $K_1$ and $K_2$
are equal. Theorem \ref{main-theorem} is proved.
\hspace{\fill}$\square$

\medskip
To illustrate the proof, let $D_1$ be the chord diagram from our
standard example. Pick a singular knot representing $D_1$, say
$$K_1 =\ \risS{-25}{kn-sin1}{}{80}{15}{8}\hspace{3cm}
  D_1 =\ \risS{-25}{cd-st-ex}{}{55}{30}{30}
$$
To perform a mutation in the share containing the chords 1,5,6, we
must slide the double point 1 close to the double points 5 and 6,
and then shrink the corresponding arcs:
$$\risS{-25}{kn-sin2}{
    \put(-5,-10){\mbox{\scriptsize Sliding the double point 1}}}{80}{20}{30}\hspace{2cm}
  \risS{-25}{kn-sin3}{
    \put(5,-10){\mbox{\scriptsize Shrinking the arcs}}}{80}{20}{30}\hspace{2cm}
  \risS{-25}{kn-sin4}{
    \put(5,-10){\mbox{\scriptsize Forming the tangle $T_S$}}
    \put(82,30){\mbox{$T_S$}}}{80}{20}{40}
$$
Now doing an appropriate rotation of the tangle $T_S$ we obtain a
singular knot $K_2$ representing the chord diagram $D_2$.

\section{Lie algebra weight systems\\
 and intersection graphs}\label{s:Lie}

Kontsevich~\cite{K} generalized a construction of
Bar-Natan~\cite{BN} of weight systems defined by a Lie algebra and
its representation to a universal weight system, with values in
the universal enveloping algebra of the Lie algebra. In~\cite{V},
Vaintrob extended this construction to Lie superalgebras.

Our main goal in this section is to prove

\begin{theorem}\label{thLie}
The universal weight systems associated to the Lie algebra
$\sL(2)$ and to the Lie superalgebra $\gl(1|1)$ depend on the
intersection graphs of chord diagrams rather than on the diagrams
themselves.
\end{theorem}

It follows immediately that the canonical knot invariants
corresponding to these two algebras do not distinguish mutants.
The latter fact is already known, but we did not manage to find
appropriate references; instead, we give a direct proof on the
intersection graphs side.

Note that for more complicated Lie algebras the statement of
Theorem~\ref{thLie} is no longer true. For example, the universal
$\sL(3)$ weight system distinguishes between the Conway and the
Kinoshita--Terasaka knots.

In fact, for each of the two algebras we prove more subtle
statements.

\begin{theorem}\label{thmatroid}
The universal weight system associated to the Lie algebra $\sL(2)$
depends on the matroid of the intersection graph of a chord
diagram rather than on the intersection graph itself.
\end{theorem}

This theorem inevitably leads to numerous questions concerning
relationship between weight systems and matroid theory, which
specialists in this theory may find worth being investigated.

Weight systems have a graph counterpart, so-called $4$-invariants
of graphs~\cite{L}. The knowledge that a weight system depends
only on the intersection graphs does not guarantee, however, that
it arises from a $4$-invariant. In particular, we do not know,
whether this is true for the universal $\sL(2)$ weight system.
Either positive (with an explicit description) or negative answer
to this question would be extremely interesting. For $\gl(1|1)$,
the answer is positive.

\begin{theorem}\label{thfour}
The universal weight systems associated to the Lie superalgebra
$\gl(1|1)$ is induced by a $4$-invariant of graphs.
\end{theorem}

In the first two subsections below, we recall the construction of
universal weight systems associated to Lie algebras and the notion
of $4$-invariant of graphs. The next two subsections are devoted to
separate treating of the Lie algebra~$\sL(2)$ and the Lie
superalgebra~$\gl(1|1)$ universal weight systems.

\subsection{Weight systems via Lie algebras}
Our approach follows that of Kontsevich in~\cite{K}. In order to
construct a weight system, we need a complex Lie algebra endowed
with a nondegenerate invariant bilinear form $(\cdot,\cdot)$. The
invariance requirement means that $(x,[y,z])=([x,y],z)$ for any
three elements~$x,y,z$ in the Lie algebra. Pick an orthonormal
basis $a_1,\dots,a_d$, $(a_i,a_j)=\delta_{ij}$, $d$ being the
dimension of the Lie algebra. Any chord diagram can be made into
an arc diagram by cutting the circle at some point and further
straightening it. For an arc diagram of~$n$ arcs, write on each
arc an index~$i$ between~$1$ and~$d$, and then write on both ends
of the arc the letter~$a_i$. Reading all the letters left to right
we obtain a word of length~$2n$ in the alphabet~$a_1,\dots,a_d$,
which is an element of the universal enveloping algebra of our Lie
algebra. The sum of all these words over all possible settings of
the indices is the element of the universal enveloping algebra
assigned to the chord diagram. This element is independent of the
choice of the cutting point of the circle, as well as the
orthonormal basis. It belongs to the center of the universal
enveloping algebra and satisfies the $4$-term relation, whence can
be extended to a weight system. The latter is called the {\it
universal weight system} associated to the Lie algebra and the
bilinear form, and it can be specialized to specific
representations of the Lie algebra as in the original Bar-Natan's
approach. Obviously, any universal weight system is
multiplicative: its value on a product of chord diagrams coincides
with the product of its values on the factors.

The simplest noncommutative Lie algebra with a nondegenerate
invariant bilinear form is $\sL(2)$. It is $3$-dimensional, and
the center of its universal enveloping algebra is the ring $\C[c]$
of polynomials in a single variable~$c$, the Casimir element. The
corresponding universal weight system was studied in detail
in~\cite{CV}. It attracts a lot of interest because of its
equivalence to the colored Jones polynomials.

In~\cite{V}, Kontsevich's construction was generalized to Lie
superalgebras, and this construction was elaborated in~\cite{FKV}
for the simplest non-commutative Lie superalgebra $\gl(1|1)$. The
center of the universal enveloping algebra of this algebra is the
ring of polynomials $\C[c,y]$ in two variables. The value of the
corresponding universal weight system on a chord diagram with~$n$
chords is a quasihomogeneous polynomial in~$c$ and~$y$, of
degree~$n$, where the weight of~$c$ is set to be~$1$, and the
weight of~$y$ is set to be~$2$.

\subsection{The $4$-bialgebra of graphs}

By a graph, we mean a finite undirected graph without loops and
multiple edges. Let $\cG_n$ denote the vector space freely spanned
over~$\C$ by all graphs with~$n$ vertices, $\cG_0=\C$ being
spanned by the empty graph. The direct sum
$$
\cG=\cG_0\oplus\cG_1\oplus\cG_2\oplus\dots
$$
carries a natural structure of a commutative cocommutative graded
Hopf algebra. The multiplication in this Hopf algebra is induced
by the disjoint union of graphs, and the comultiplication is
induced by the operation taking a graph~$G$ into the sum $\sum
G_U\otimes G_{\bar{U}}$, where $U$ is an arbitrary subset of
vertices of~$G$, $\bar{U}$ its complement, and $G_U$ denotes the
subgraph of~$G$ induced by~$U$.

The {\it $4$-term relation for graphs} is defined in the following
way. By definition, the {\it $4$-term element} in~$\cG_n$
determined by a graph~$G$ with~$n$ vertices and an ordered pair
$A,B$ of its vertices connected by an edge is the linear combination
$$
G-G'_{AB}-\widetilde{G}_{AB}+\widetilde{G}'_{AB},
$$
where
\begin{itemize}
\item $G'_{AB}$ is the graph obtained by deleting the edge~$AB$
in~$G$;

\item $\widetilde{G}_{AB}$ is the graph obtained by switching the
adjacency to~$A$ of all the vertices adjacent to~$B$ in~$G$;

\item $\widetilde{G}'_{AB}$ is the graph obtained by deleting the
edge~$AB$ in~$G'_{AB}$ (or, equivalently, by switching the
adjacency to~$A$ of all the vertices adjacent to~$B$
in~$G'_{AB}$).
\end{itemize}

All the four terms in a~4-term element have the same number~$n$ of
vertices. The quotient of~$\cG_n$ modulo the span of all 4-term
elements in~$\cG_n$ (defined by all graphs and all ordered pairs
of adjacent vertices in each graph) is denoted by~$\cF_n$. The
direct sum
$$
\cF=\cF_0\oplus\cF_1\oplus\cF_2\oplus\dots
$$
is the quotient Hopf algebra of graphs, called the {\it
$4$-bialgebra}. The mapping taking a chord diagram to its
intersection graph extends to a graded Hopf algebra homomorphism
$\gamma$ from the Hopf algebra of chord diagrams to $\cF$.

Being commutative and cocommutative, the $4$-bialgebra is
isomorphic to the polynomial ring in its basic primitive elements,
that is, it is the tensor product $S(\cP_1)\otimes
S(\cP_2)\otimes\dots$ of the symmetric algebras of its homogeneous
primitive spaces.

\subsection{The $\sL(2)$ weight system}\label{sl2}

Our treatment of the universal weight system associated with the
Lie algebra $\sL(2)$ is based on the recurrence formula for
computing the value of this weight system on chord diagrams due to
Chmutov and Varchenko~\cite{CV}. The recurrence states that if a
chord diagram contains a leaf, that is, a chord intersecting only
one other chord, then the value of the $\sL(2)$ universal weight
system on the diagram is~$(c-1/2)$ times its value on the result of
deleting the leaf, and, in addition,
\def\chvad#1{\risS{-10}{#1}{}{25}{10}{10}}
$$\chvad{chvard10}\ -\ \chvad{chvard11}\ -\
     \chvad{chvard12}\ +\  \chvad{chvard13}\ =\
  2\ \chvad{chvard14}\ -\ 2\ \chvad{chvard15}
$$
meaning that the value of the weight system on the chord
diagram on the left-hand side coincides with the linear
combinations of its values on the chord diagrams indicated on the
right.

Now, in order to prove Theorem~\ref{thLie} for the universal
$\sL(2)$ weight system, we must prove that mutations of a chord
diagram preserve the values of this weight system. Take a chord
diagram and a share in it. Apply the above reccurence formula to a
chord and two its neighbors belonging to the chosen share. The
recurrence relation does not affect the complementary share, while
all the instances of the modified first share are simpler than the
initial one (each of them contains either fewer chords or the same
number of chords but with fewer intersections). Repeating this
process, we can replace the original share by a linear combination
of the simplest shares, chains, which are symmetric meaning that
they remain unchanged under rotations. The $\sL(2)$ case of
Theorem~\ref{thLie} is proved.
\hspace{\fill}$\square$

\bigskip
Now let us turn to the proof of Theorem~\ref{thmatroid}. For
elementary notions of matroid theory we refer the reader to any
standard reference, say to~\cite{W}. Recall that a matroid can be
associated to any graph. It is easy to check that the matroid
associated to the disjoint union of two graphs coincides with that
for the graph obtained by identifying a vertex in the first graph
with a vertex in the second one. We call the result of gluing a
vertex in a graph~$G_1$ to a vertex in a graph~$G_2$ a {\it
$1$-product\/} of $G_1$ and~$G_2$. The converse operation is {\it
$1$-deletion}. Of course, the $1$-product depends on the choice of
the vertices in each of the factors, but the corresponding matroid
is independent of this choice.

Similarly, let $G_1$, $G_2$ be two graphs, and pick vertices
$u_1,v_1$ in~$G_1$ and $u_2,v_2$ in~$G_2$. Then the matroid
associated to the graph obtained by identifying $u_1$ with $u_2$
and $v_1$ with $v_2$ coincides with the one associated to the
graph obtained by identifying $u_1$ with $v_2$ and $u_2$ with
$v_1$. The operation taking the result of the first identification
to that of the second one is called the {\it Whitney twist\/} on
graphs.

Both the $1$-product and the Whitney twist have chord diagram
analogs. For two chord diagrams with a distinguished chord in each
of them, we define their $1$-product as a chord diagram obtained
by replacing the distinguished chords in the ordinary product of
chord diagrams chosen so as to make them neighbors by a single
chord connecting their other ends. The Whitney twist also is well
defined because of the following statement.

\begin{lemma}
Suppose the intersection graph of a chord diagram is the result of
identifying two pairs of vertices in two graphs~$G_1$ and~$G_2$.
Then both graphs~$G_1$ and~$G_2$ are intersection graphs, as well
as the Whitney twist of the original graph.
\end{lemma}

The assertion concerning the graphs~$G_1$ and~$G_2$ is obvious. In
order to prove that the result of the Whitney twist also is an
intersection graph, let $c_1,c_2$ denote the two chords in a chord
diagram~$C$ such that deleting these chords makes~$C$ into an
ordinary product of two chord diagrams~$C_1,C_2$. By reflecting
the diagram~$C_2$ and restoring the chords~$c_1$ and~$c_2$ we
obtain a chord diagram whose intersection graph is the result of
the desired Whitney twist. The lemma is proved.

According to the Whitney theorem, two graphs have the same matroid
iff they can be obtained from one another by a sequence of
$1$-products/deletions and Whitney twists. Therefore,
Theorem~\ref{thmatroid} follows from

\begin{lemma} {\rm (i)} The value of the universal
$\sL(2)$ weight system on the $1$-product of chord diagrams
coincides with the product of its values on the factors divided
by~$c$. {\rm (ii)} The value of the universal $\sL(2)$ weight
system remains unchanged under the Whitney twist of the chord
diagram.
\end{lemma}

Statement~(i) is proved in~\cite{CV}. The proof of statement~(ii)
is similar to that of Theorem~\ref{thLie}. Consider the part~$C_2$
participating in the Whitney twist and apply to it the recurrence
relations. Note that the relations do not affect the complementary
diagram~$C_1$. Simplifying the part~$C_2$ we reduce it to a linear
combination of the simplest possible diagrams, chains, which are
symmetric under reflection. Reflecting a chain preserves the chord
diagram, whence the value of the $\sL(2)$ weight system.
Theorem~\ref{thmatroid} is proved.
\hspace{\fill}$\square$

\subsection{The $\gl(1|1)$ weight system}\label{gl11}

Define the (unframed) {\it Conway graph invariant} with values in the ring
of polynomials~$\C[y]$ in one variable~$y$ in the following way.
We set it equal to $(-y)^{n/2}$ on graphs with $n$ vertices if the
adjacency matrix of the graph is nondegenerate, and $0$ otherwise.
Recall that the adjacency matrix~$A_G$ of a graph~$G$ with~$n$
vertices is an $n\times n$-matrix with entries in $\Z_2$ obtained
as follows. We choose an arbitrary numbering of the vertices of
the graph, and the entry $a_{ij}$ is~$1$ provided the $i$~th and
the $j$~th vertices are adjacent and~$0$ otherwise (diagonal
elements $a_{ii}$ are~$0$). Note that for odd~$n$, the adjacency
matrix cannot be nondegenerate, hence the values indeed are in the
ring of polynomials. The Conway graph invariant is
multiplicative: its value on the disjoint union of graphs is the
product of its values on the factors.

Clearly, the Conway graph invariant is a $4$-invariant. Moreover,
it satisfies the $2$-term relation, which is more restrictive than
the $4$-term one: its values on the graphs~$G$ and
$\widetilde{G}_{AB}$ coincide for any graph~$G$ and any pair of
ordered vertices~$A,B$ in it. Indeed, consider the graph as a
symmetric bilinear form on the $\Z_2$-vector space whose basis is
the set of vertices of the graph, the adjacency matrix being the
matrix of the bilinear form in this basis. In these terms, the
transformation $G\mapsto\widetilde{G}_{AB}$ preserves the vector
space and the bilinear form, but changes the basis $A,B,C,\dots\to
A+B,B,C,\dots$. Thus, it preserves the nondegeneracy property of
the adjacency matrix.

The subspace~$\cF_1$ is spanned by the graph~$p_1$ with a single
vertex (whence no edges), which is a primitive element.
Since~$\cF$ is the polynomial ring in its primitive elements, each
homogeneous space~$\cF_n$ admits a decomposition into the direct
sum of two subspaces, one of which is the subspace of polynomials
in primitive elements of degree greater than~$1$, and the other
one is the space of polynomials divisible by~$p_1$. We define the
{\it framed Conway graph invariant} as the only multiplicative
$4$-invariant with values in the polynomial ring $\C[c,y]$ whose
value on~$p_1$ is~$c$, and on the projection of any graph to the
subspace of $p_1$-independent polynomials along the subspace of
$p_1$-divisible polynomials coincides with the Conway graph
invariant of the graph.

The values of the framed Conway graph invariant can be computed
recursively. Take a graph~$G$ and consider its projection to the
subspace of graphs divisible by~$p_1$. On this projection, the
framed Conway graph invariant can be computed because of its
multiplicativity. Now add to the result the value of the
(unframed) Conway graph invariant on the graph. Now we can refine the 
statement of theorem \ref{thfour}.

\begin{theorem}
The $\gl(1|1)$ universal weight system is the pullback of the
framed Conway graph invariant to chord diagrams under the
homomorphism~$\gamma$.
\end{theorem}

{\bf Proof.} The proof follows from two statements in~\cite{FKV}.
Theorem~3.6 there states that setting $c=0$ in the value of the
$\gl(1|1)$ universal weight system on a chord diagram we obtain
the result of deframing this weight system. Theorem~4.4 asserts
that this value is exactly the Conway invariant of the chord
diagram. The latter coincides with the Conway graph invariant of
the intersection graph of the chord diagrams defined above. Since
the deframing for chord diagrams is a pullback of the deframing
for graphs, we are done.
\hspace{\fill}$\square$

\end{document}